\documentclass[runningheads]{llncs}
\usepackage[utf8]{inputenc}
\usepackage{graphicx}
\usepackage{lmodern}
\usepackage{
    multirow,     
    multicol,     
    amsmath,      
    amsfonts,     
    amssymb,      
    bm,           
    siunitx,      
    wasysym,      
    float,        
    enumitem,     %
    eurosym,      %
    pgfplotstable,%
    makecell,
    booktabs,
    colortbl,
    algorithm,
    algorithmic,
    geometry,
    tikz      %
}
\usetikzlibrary{calc}
\usepackage[colorlinks=true,allcolors=blue]{hyperref}
\usepackage{csquotes}

\usepackage[english]{babel}
\SetLabelAlign{parright}{\parbox[t]{\labelwidth}{\raggedleft#1}}

\setlist[description]{style=multiline,topsep=10pt,leftmargin=2cm,font=\normalfont,align=parright}

\usepackage{svg}
\usepackage[long,,nocomma]{optidef}

\usepackage[super]{nth}

\newcommand{\ProbOpt}[1]
{
	$\mathcal{P}_\text{#1}$
}

\usepackage{epstopdf}
\usepackage{subfigure}

\begin{document}
\title{An Application of BnB-NSGAII: Initializing NSGAII to Solve 3 Stage Reducer Problem  \thanks{Supported by organization ERDF, Grand Est and Lebanese University}}
\titlerunning{An Application of BnB-NSGAII}
%
\author{Ahmed Jaber\inst{1,2}\orcidID{0000-0001-5508-1299} \and
Pascal Lafon\inst{1} \and
Rafic Younes\inst{2}}
\authorrunning{A. Jaber et al.}
%
\institute{University of Technology of Troyes UTT, France\\
\url{https://www.utt.fr/} \and
Lebanese University, Beirut, Lebanon\\
\url{https://ul.edu.lb/}}
\maketitle              

\begin{abstract}
The 3 stage reducer problem is a point of interest for many researchers. In this paper, this problem is reformulated to a bi-objective problem with additional constraints to meet the ISO mechanical standards. Those additional constraints increase the complexity of the problem, such that, NSGAII performance is not sufficient. To overcome this, we propose to use BnB-NSGAII \cite{jaber_branch-and-bound_2021} method - a hybrid multi-criteria branch and bound with NSGAII - to initialize NSGAII before solving the problem, seeking for a better initial population. A new feature is also proposed to enhance BnB-NSGAII method, called the legacy feature. The legacy feature permits the inheritance of the elite individuals between - branch and bound - parent and children nodes. NSGAII and BnB-NSGAII with and without the legacy feature are tested on the 3 stage reducer problem. Results demonstrate the competitive performance of BnB-NSGAII with the legacy feature.
\end{abstract}
\keywords{NSGAII  \and multi-objective \and MINLP \and branch-and-bound \and 3-stage reducer.}

\section{Introduction}
 In \cite{Fauroux1999}, the design of the 3 Stage Reducer (3SR) optimization problem has been introduced to illustrate the optimal design framework of the power transmission mechanism. This problem has been a point of interest for many researchers in different domains. Engineering researchers enhance the problem for mechanical engineering applications. In  \cite{Fauroux2004}, the problem is extended to a mixed variables optimization problem. And recently a similar problem is stated in \cite{Han2020} to illustrate the optimization of the volume and layout design of 3SR. Due to the problem complexity, optimization researchers are interested to test optimization methods on it. In \cite{Yvars2018}, the authors use the 3SR problem to examine the performance of the constraint propagation method.\par
 In this paper, the 3SR problem is reformulated to a bi-objective problem with additional constraints to meet the ISO mechanical standards. Those additional constraints increase the complexity of the problem, such that, the well-known Non-Dominated Sorting Genetic Algorithm 2 (NSGAII) \cite{deb_fast_2002} performance is not sufficient.\par
 In \cite{jaber_branch-and-bound_2021}, the authors enhance the performance of NSGAII by hybridizing it with the multi-criteria branch and bound method \cite{mavrotas_branch_1998}, the proposed method is called BnB-NSGAII. In this paper, we propose to use the BnB-NSGAII method to initialize NSGAII before solving the 3SR problem, seeking a better initial population. The initial population seeding phase is the first phase of any genetic algorithm application. It generates a set of solutions randomly or by heuristic initialization as input for the algorithm. Although the initial population seeding phase is executed only once, it has an important role to improve the genetic algorithm performance \cite{deng_improved_2015}.\par 
 Furthermore, we propose a new feature to enhance the BnB-NSGAII method, called the legacy feature. The legacy feature permits the inheritance of elite genes between branch-and-bound nodes.\par
 The rest is organized as follows. Section \ref{3st} presents the 3SR problem and its complexity. The proposed BnB-NSGAII legacy feature is explained in section \ref{bnb_nsgaii}. The computational results are reported in section \ref{test}. Finally, an overall conclusion is drawn in section \ref{conc}.
\section{3 Stage Reducer Problem} \label{3st}
\begin{figure}[h!]
\centering
\includegraphics[width=0.55\linewidth]{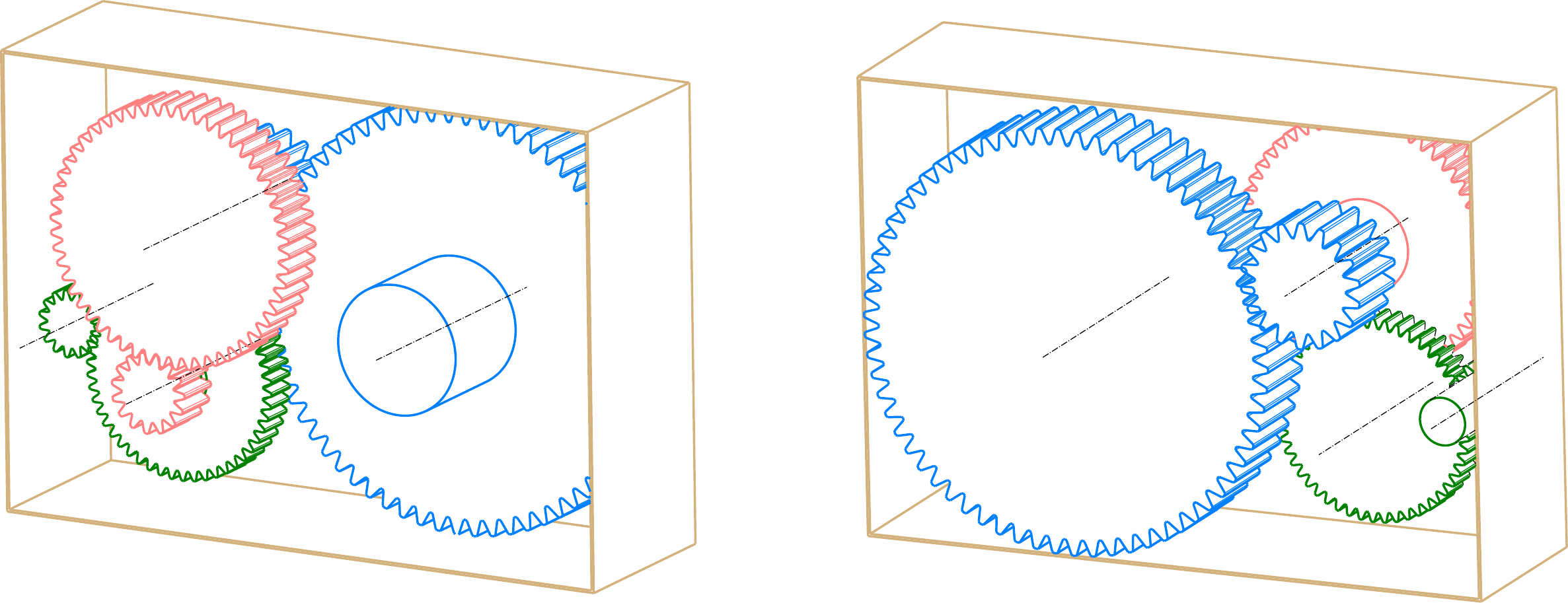}
\caption{Front and back view of a 3 stage reducer with closure.}
\label{fig:reducteuriftomm}
\end{figure}

\begin{figure}[h!]
\centering
\includegraphics[width=0.45\linewidth]{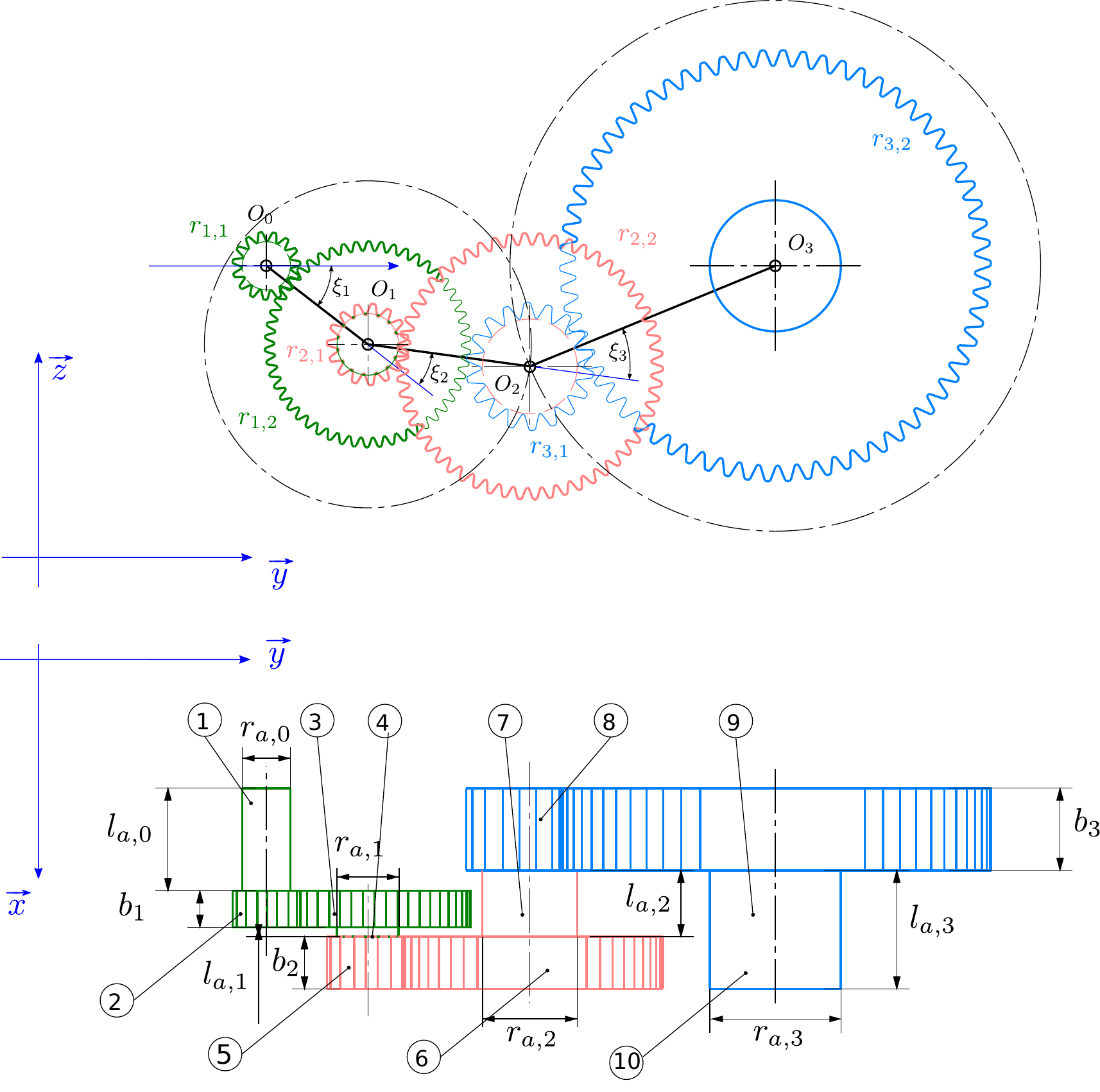}
\caption{Detailed view of the 3 stage reducer.}
\label{fig:reducteuriftommvueface}
\end{figure}

 The design problem consists in finding dimensions of  main components (pinions, wheels and shafts) of the 3 stage reducer (figure \ref{fig:reducteuriftomm}) to minimize the following bi-objective problem :
\begin{enumerate}
\item The volume of all the components of the reducer :
\begin{equation}
f_1(\bm{x})=\pi\left(\sum\limits_{s=0}^{s=3}l_{a_s}{r_{a,s}}^2 
+ 
\sum\limits_{s=1}^{s=3}\left[
b_s\frac{m_{\text{n}s}^2}{2}(Z_{s,1}^2+Z_{s,2}^2)
\right]
\right)
\end{equation}
\item The gap between the required reduction ratio $\bar{u}$ and the ratio of the reducer (tolerance):
\begin{equation}
f_2(\bm{x}) =\frac{1}{\bar{u}} \left\vert\bar{u}-\prod\limits_{s=1}^{s=3}\frac{Z_{s,2}}{Z_{s,1}}\right\vert, \;\bar{u}>1
\end{equation}
\end{enumerate}
The problem is designed assuming the following are known:
    \begin{itemize}
    \item The power to be transmitted, $P_t$ and the speed rotation of input shaft $N_e$.
    \item The total speed rotation reduction ratio $\bar{u}$, the position of the output shaft from the input shaft position (figure \ref{fig:reducteuriftommvueface}).
    \item The dimension of the casing box.
    \end{itemize}

The 3SR problem is formulated with 2 objective functions, 41 constraints (presented in Appendix \ref{prob_cnt}), 3 categorical variables (gears modules), 6 integer variables (number of teeth), and 11 continuous variables. Gears modules have  41 possibilities, pinion number of teeth ranges from 14 to 30 and wheel number of teeth ranges from 14 to 150. Hence, the combinatorial space of the 3SR problem consists in $41^3+(30-14)^3+(150-14)^3 \simeq \SI{8.7E14}{} $. Thus, the problem is considered a mid-sized problem concerning the number of variables and constraints, but, huge combinatorial space.\par
\begin{figure}[h!]
\centering
\subfigure[Number of converged runs.]{\includegraphics[width=0.48\linewidth]{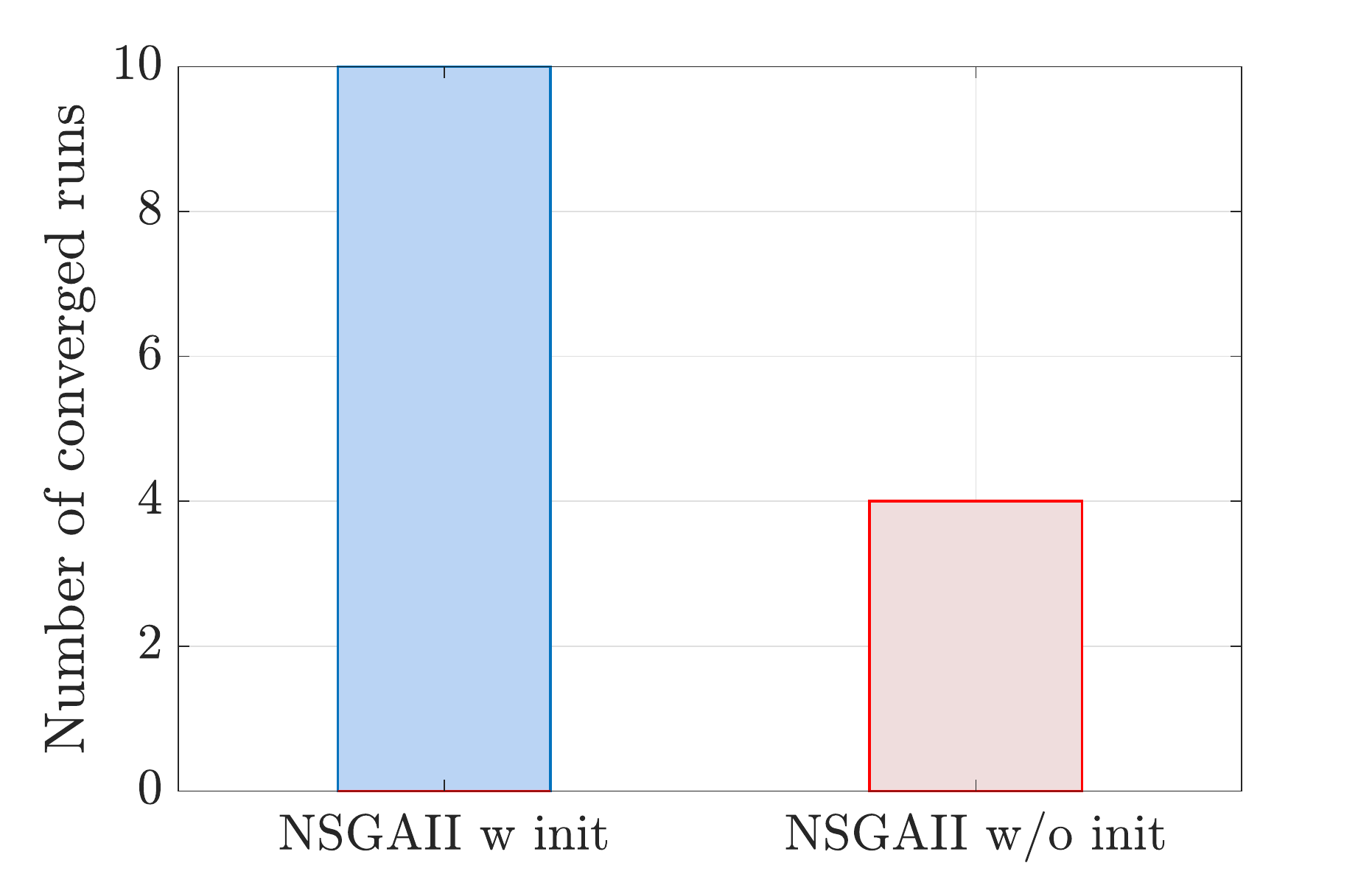}
    \label{fig:param} }
\subfigure[Best Pareto front.]{\includegraphics[width=0.48\linewidth]{ 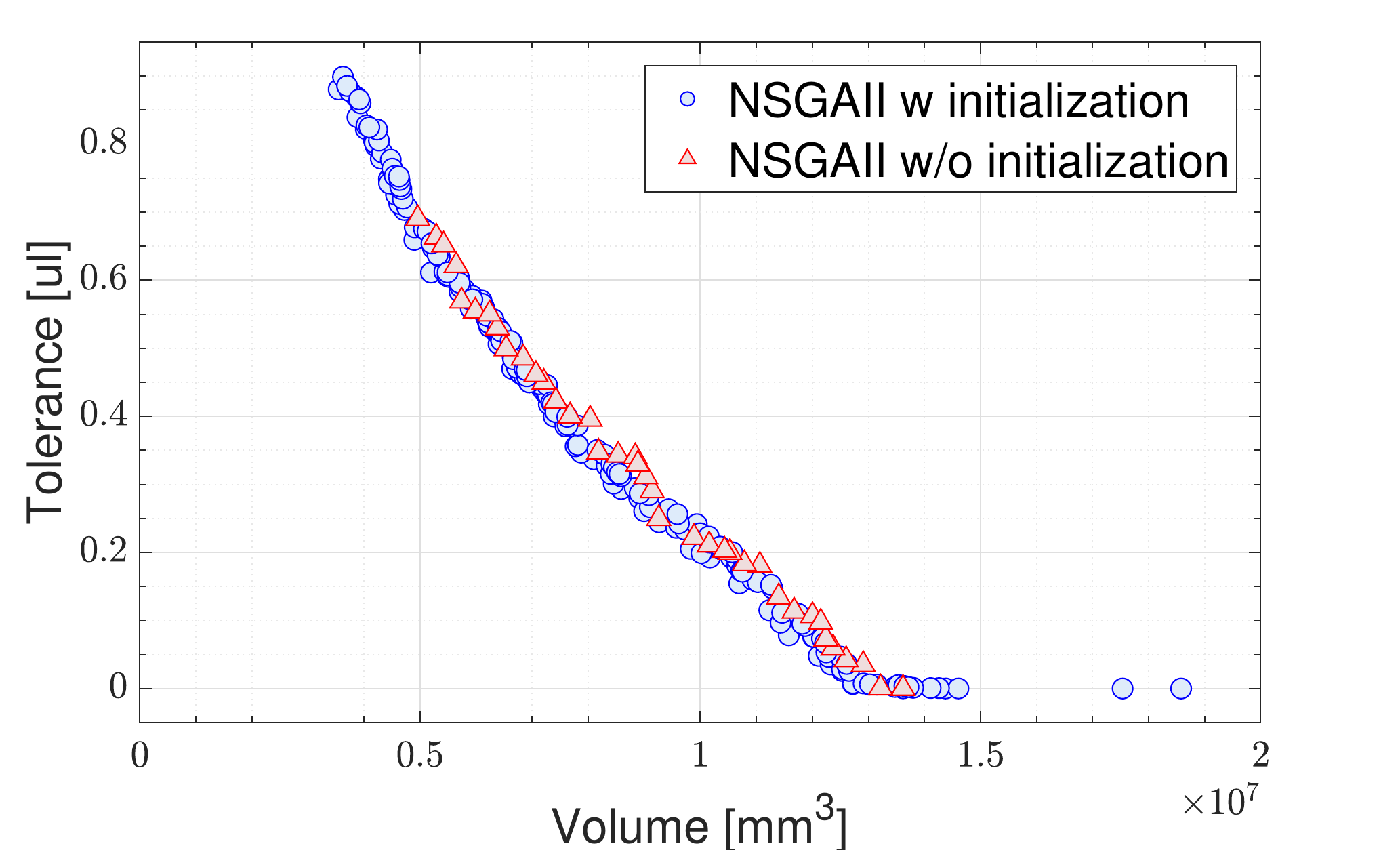}
\label{fig:nsga_init}}
\caption{Results of 3SR problem solved by NSGAII with (blue) and without (red) initial feasible seed.}
\end{figure}

The additional constraints increase the complexity of the problem. This is noticed by solving the problem using NSGAII with different initial conditions as follows. In first hand, NSGAII is initialized with 1 feasible individual. On the other hand, NSGAII is randomly initialized. Each was run 10 times with the same parameters shown in Table \ref{table_param}. Figure \ref{fig:param} shows how many run each method converged to a feasible solution out of 10. Figures \ref{fig:param} and \ref{fig:nsga_init} show that if the initial population contains at least 1 feasible individual, NSGAII converges to a good approximated Pareto front every time. Whilst, if NSGAII is initialized with a random population, NSGAII either fails to converge to a feasible solution, or it converges to a low-quality Pareto front.\par
\begin{table}[H]
\centering
	\caption{Parameters used for NSGAII algorithm}
	\label{table_param}
	\begin{tabular}{ll}
		\hline
		\textbf{Parameters}               & \multicolumn{1}{l}{\textbf{Value}}  \\ \hline
		Cross over probability             & 0.8                                \\
		Mutation Probability              & 0.9                                \\
		Population size                   & 200               \\
		Allowable generations     & 500             \\
		
		Constraint handling     & Legacy method \cite{deb_fast_2002}                     \\ 
		Crossover operator				& Simulated Binary crossover (SBX)  \cite{maruyama_parametric_2017}   \\
		ETAC & 100
		\\
		 Mutation operator				&  Partially-mapped crossover (PMX) \cite{maruyama_parametric_2017} \\
		ETAM & 10
		\\
		 \hline
	\end{tabular}
\end{table}

\begin{figure}[hb]
\centering
\includegraphics[width=0.8\linewidth]{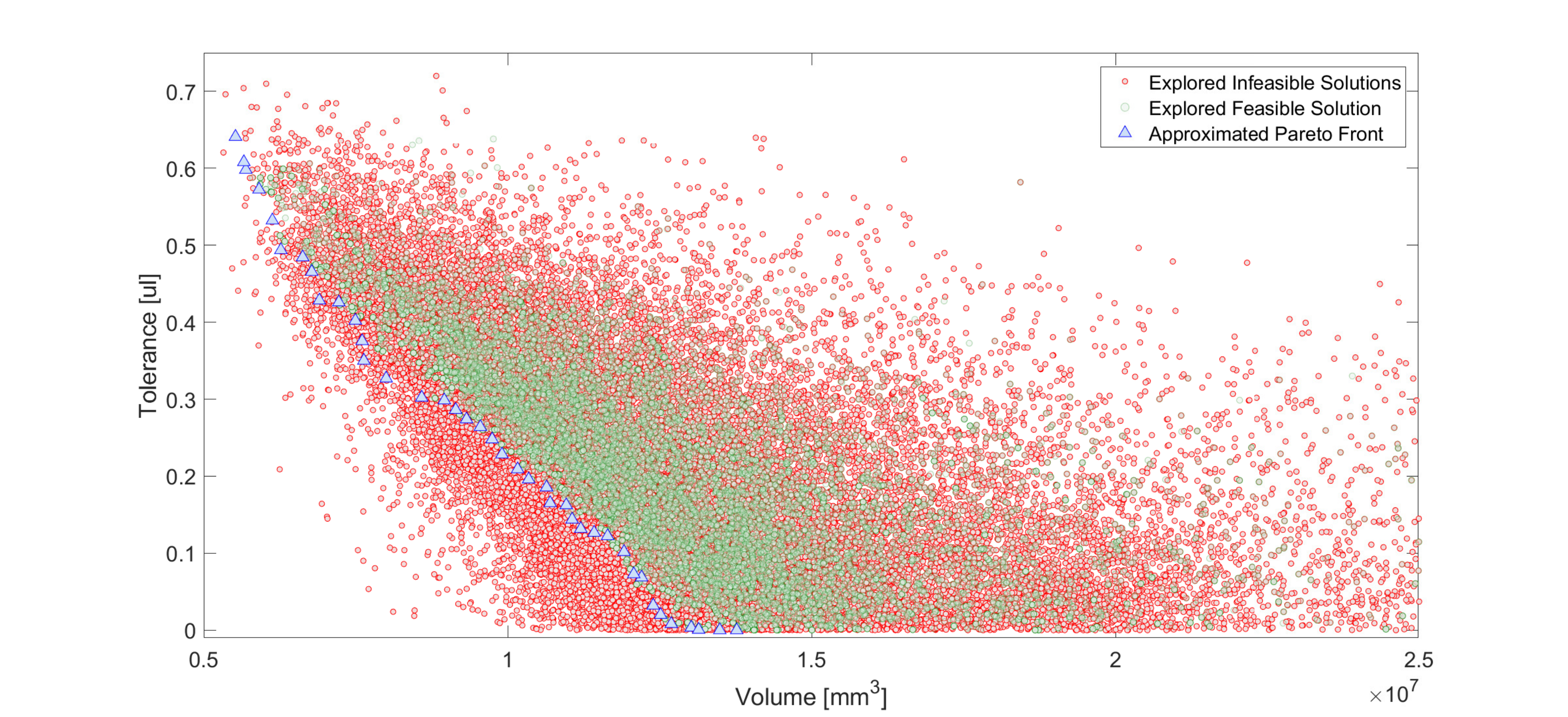}
\caption{Explored portion of the domain, showing the 3SR problem complexity.}
\label{fig:prob_cmplx}
\end{figure}

Figure \ref{fig:prob_cmplx} shows part of the domain of the 3SR problem explored by NSGAII with feasible initial population. The explored domain shows the complexity of the problem, where both feasible and infeasible solutions share the same domain on the projected objective domain. Moreover, all the feasible solutions are too near to the infeasible ones.\par
To enhance the quality of the solution of this problem - and accordingly any similar problem - where feasible solutions are not known, our proposal is first to use BnB-NSGAII proposed in \cite{jaber_branch-and-bound_2021} to search for feasible individuals. These individuals are then injected in the random initial population of NSGAII.\par
\section{BnB-NSGAII} \label{bnb_nsgaii}
In \cite{jaber_branch-and-bound_2021} the authors proposed the BnB-NSGAII approach. In this approach, Multi-Criteria Branch and Bound (MCBB) \cite{mavrotas_multi-criteria_2005} is used to enhance the exploration force of NSGAII by investigating the mixed-integer domain space through branching it to subdomains, then NSGAII bounds each one. In this way, MCBB guides the search using the lower bounds obtained by NSGAII. Our proposal is to furthermore enhance the exploration potential of BnB-NSGAII by adding the legacy feature.\par
\subsection{General Concept of BnB-NSGAII}
The general multi-Objective MINLP problem (\ProbOpt{MO-MINLP}) is written as
\begin{mini}
        {\bm{x},\bm{y}}{\bm{f}(\bm{x}, \bm{y}) = {f_1(\bm{x}, \bm{y}),\ldots, f_p(\bm{x}, \bm{y})}}{}{}
        \addConstraint{}
        \label{mominlp}
        \addConstraint{c_j(\bm{x}, \bm{y})\leq0,\; }{}{j=1,...,m}
        \addConstraint{\bm{x} \in \bm{X},}{}{\bm{X}\in \mathbb{R}^{n_c}}
        \addConstraint{\bm{y} \in \bm{Y},}{}{\bm{Y}\in \mathbb{N}^{n_i},}
\end{mini}
 where $p$ and $m$ are the number of objectives and constraints respectively. $\bm{X}$ and $\bm{Y}$ denote the set of feasible solutions of the problem for $n_c$ continuous and $n_i$ integer variables respectively. \par
\ProbOpt{MO-MINLP} is complex and expensive to solve. The general idea is thus to solve several simpler problems instead. BnB-NSGAII divides \ProbOpt{MO-MINLP} by constructing a combinatorial tree that aim to partition the root node problem - \ProbOpt{MO-MINLP} - into a finite number of subproblems $Pr_1,\dots,Pr_i,\dots,Pr_n$. Where $i$ and $n$ are the current node and the total number of nodes respectively. Each $Pr_i$ is considered a node. Each node is then solved by NSGAII. Solving a node is to determine its lower and upper bounds. The upper bound of a node $P^N_i$ is the Pareto front captured by NSGAII, which is then stored in an incumbent list $P^N$. Whilst the lower bound is the ideal point $P^I_i$ of the current node. 
 \begin{equation}\label{idealeq}
     P^I_i=\min f_k(\bm{x}_i,\bm{y}_i); \;\;\;\;k=1,\dots,p.
 \end{equation}
By solving $Pr_i$, one of the following is revealed:
\begin{itemize}
        \item[$\bullet$]  $Pr_i$ is infeasible, means that NSGAII didn't find any solution that satisfies all constraints. Hence, $Pr_i$ is pruned (\textit{fathomed}) by \textit{infeasibility}.
        \item[$\bullet$] $Pr_i$ is feasible, but, the current lower bound $P^I_i$ is dominated by a previously found upper bound $P^N$. Therefore, $Pr_i$ is fathomed by \textit{optimality}.
        \item[$\bullet$] $Pr_i$ is feasible, and, $P^I_i$ is not dominated by $P^N$, $P^I_i \leq P^N$. $P^N$ is then updated by adding $P^N_i$ to it.
\end{itemize}
In the \nth{3} case, the combinatorial tree is furtherly branched by dividing $Pr_i$ into farther subproblems, called children nodes. If a node cannot be divided anymore, it is called a leaf node. Leaves are the simplest nodes, since all integer variables are fixed such that $\bm{y}=\bm{\bar{y}}$. NSGAII then solve leaves as Multi-Objective continuous Non-Linear problem (\ProbOpt{MO-NLP}): 
\begin{mini}
        {\bm{x},\bm{\bar{y}}}{\bm{f}(\bm{x}, \bm{\bar{y}}) = {f_1(\bm{x},  \bm{\bar{y}}),\ldots, f_p(\bm{x}, \bm{\bar{y}})}}{}{}
        \addConstraint{}
        \label{monlp}
        \addConstraint{c_j(\bm{x}, \bm{\bar{y}})\leq0,\;}{}{j=1,...,m}
        \addConstraint{}{\bm{x} \in \bm{X}_i,}{}
\end{mini} 
where $\bm{X}_i$ denotes the set of feasible solutions of the current node. $P^N_i$ of each leaf is then added to $P^N$. The overall Pareto front is obtained by removing the dominated elements from $P^N$.
\subsection{BnB Legacy Feature}
In NSGAII, the best population is that found in the last generation, since it contains the elite individuals among all the previous generations. In BnB-NSGAII, each node is solved independently. The output of each node is the captured Pareto front only. The last population in the node is thus discarded, although it might be valuable to other nodes.\par 
We propose to permit the legacy between nodes. Where each child node inherits the last population from its parent node. The child node then initializes NSGAII by this population.\par 
The children nodes are subproblems of their parent node. Thus, the boundary of parent node is different than that for the children nodes, $\bm{Y}_{parent} \neq \bm{Y}_{child}$. Hence, the population is rebounded before initializing NSGAII. Rebounding the population may lead to the loss of the elite individuals, though some of the elite genes are still conserved.\par
\subsection{An Application of BnB-NSGAII}
\begin{figure}[h!]
    \centering
    \includegraphics[width=0.85\linewidth]{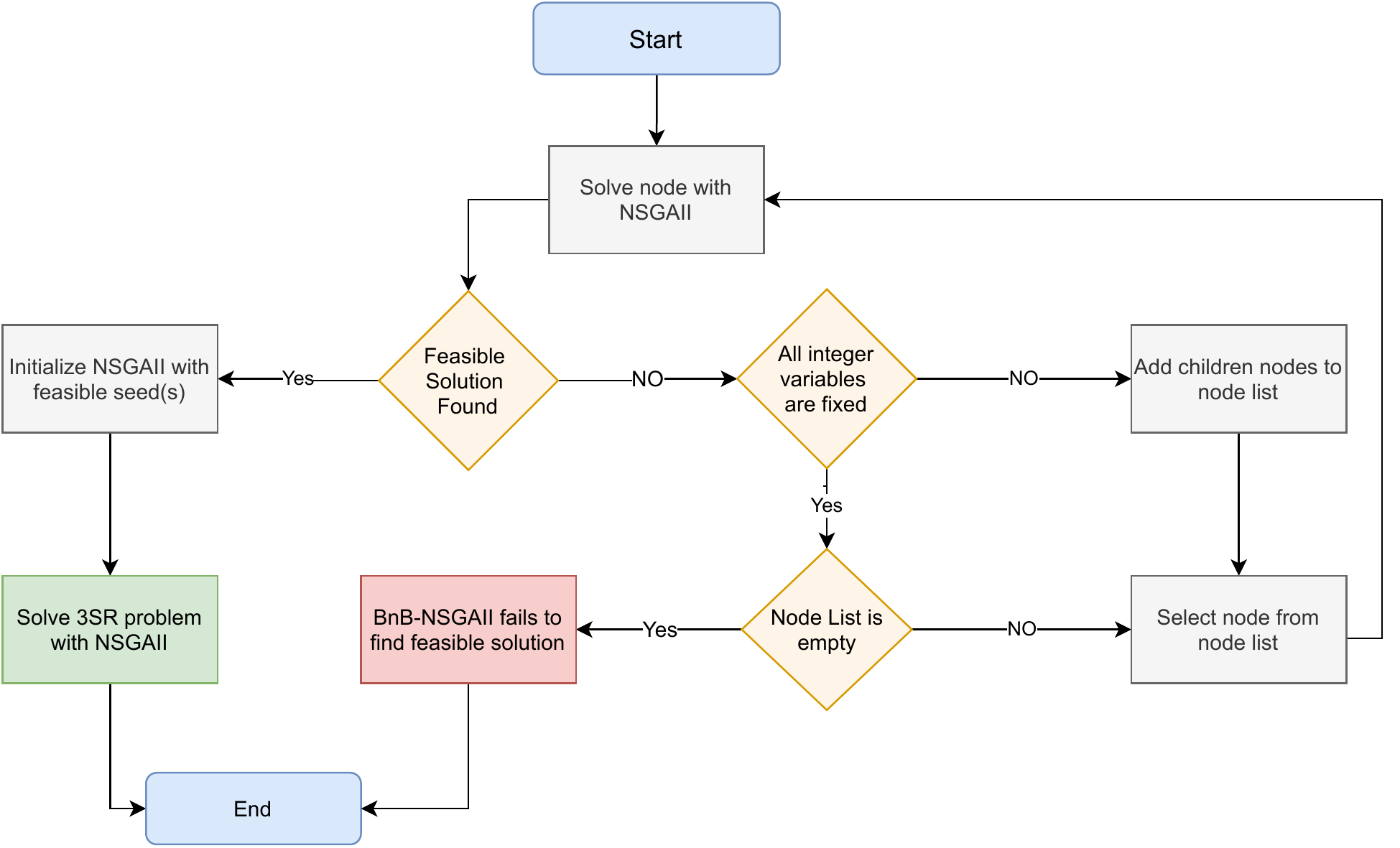}
    \caption{Flowchart of BnB-NSGAII application.}
    \label{fig:my_label}
\end{figure}

BnB-NSGAII is characterized by high exploration potential. Thus, in this paper, BnB-NSGAII is used to search for at least one feasible solution for the 3SR problem. For this aim, BnB-NSGAII is properly modified to 1) continue enumeration of the combinatorial tree even if the root node is infeasible. 2) stop whenever a feasible solution(s) is found. Then, NSGAII is called to solve the 3SR problem by initializing it with the feasible solution(s) found as shown in Figure \ref{fig:my_label}.
\section{Numerical Experiment}\label{test}
NSGAII and BnB-NSGAII with and without the legacy feature were tested on the 3SR problem. Each method was run 10 times. The test was done using the same parameters for the 3 solvers. Table \ref{table_param} shows the parameters used in this experiment.\par
\subsection{Results and Discussion}
\begin{figure}[h!]
    \centering
    \subfigure[Number of converged runs]{\includegraphics[width=0.48\linewidth, height=0.31\linewidth]{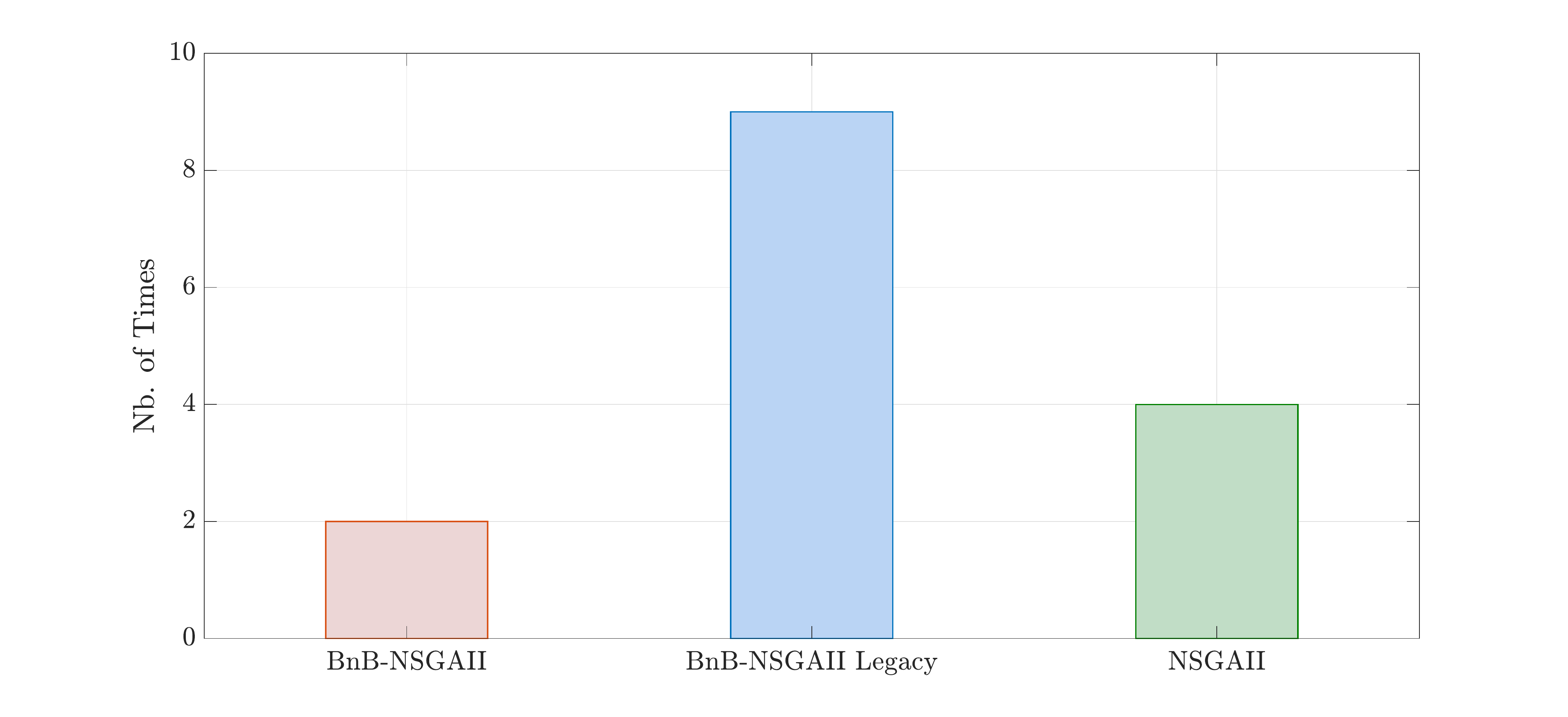}
    \label{fig:stat}
    }
    \subfigure[NSGAII]{\includegraphics[width=0.48\linewidth]{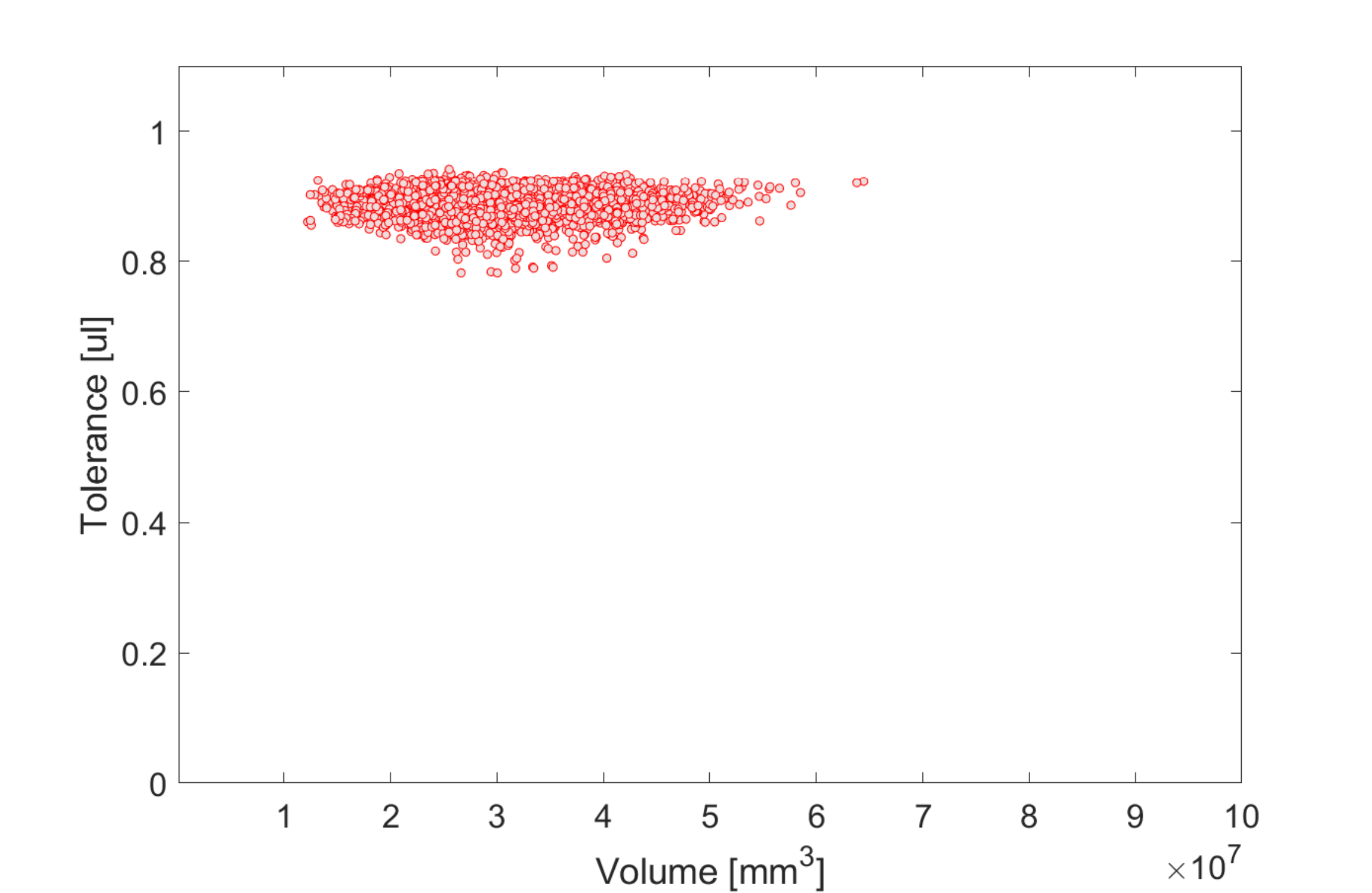}
    \label{fig:nsga_failed}
    }
    \subfigure[BnB-NSGAII]{\includegraphics[width=0.48\linewidth]{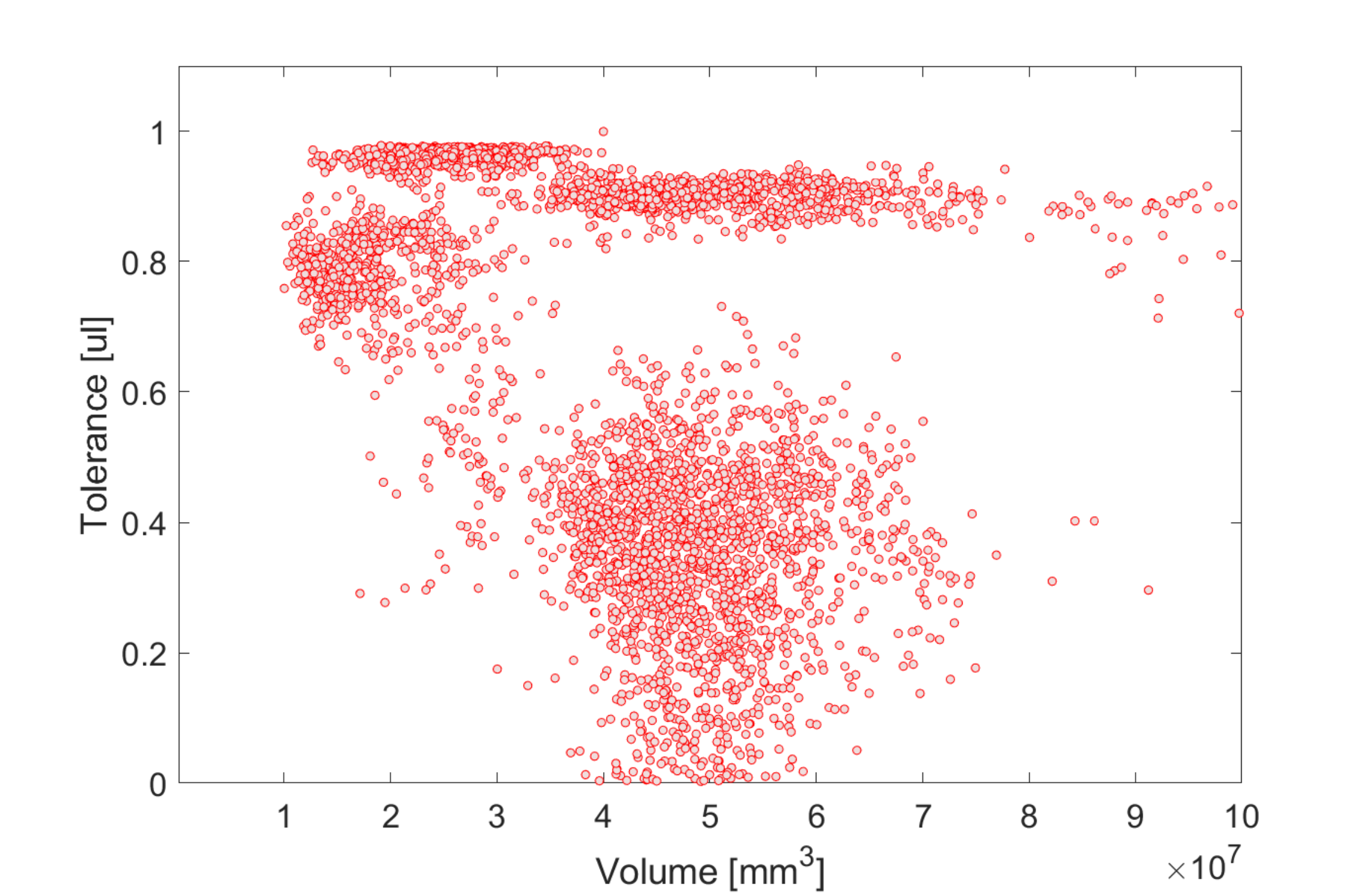}
    \label{fig:bnb_nsga_failed}
    }
    \subfigure[BnB-NSGAII Legacy]{\includegraphics[width=0.48\linewidth]{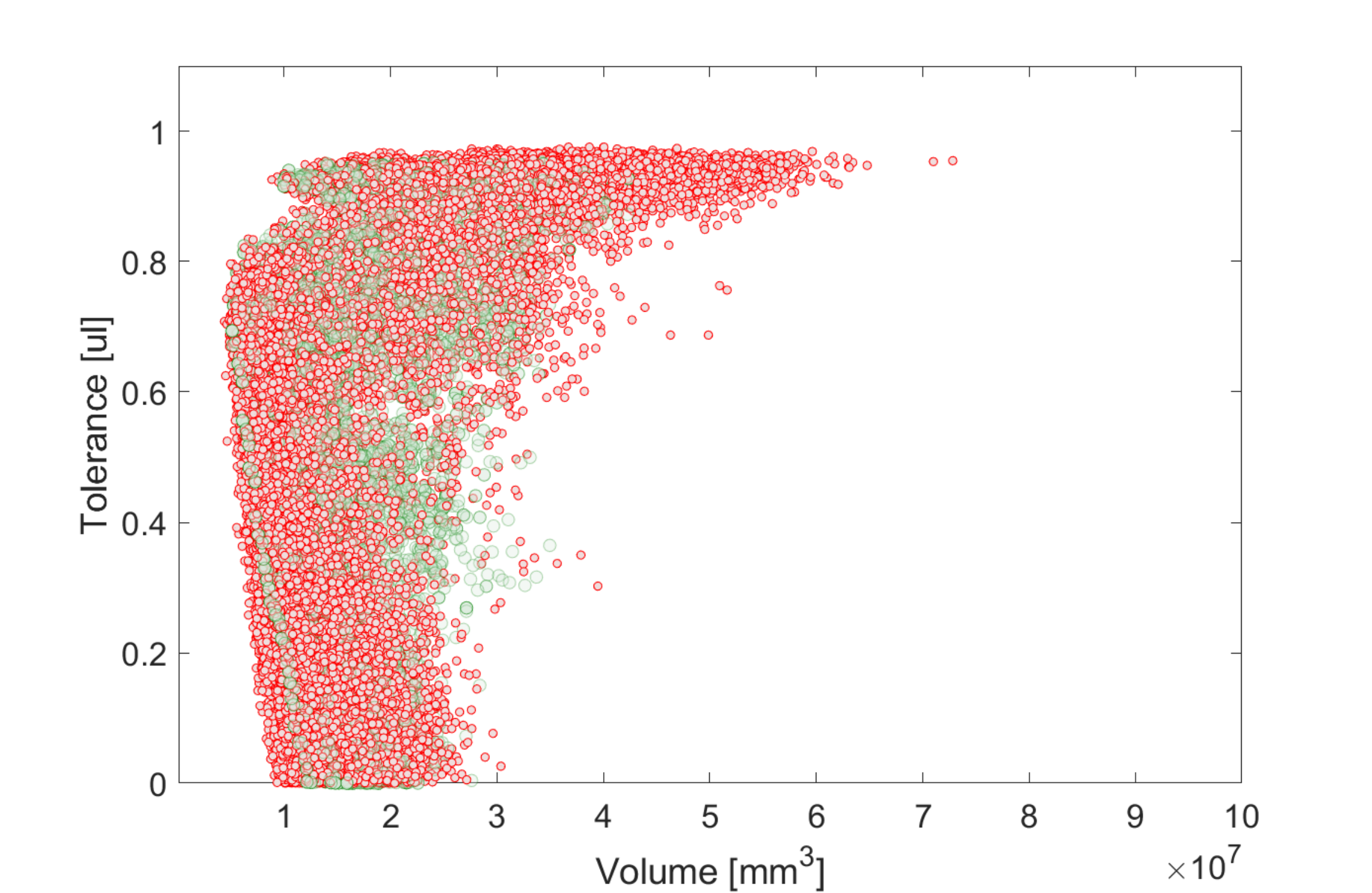}
    \label{fig:bnb_legacy}
    }
    \caption{Explored domain by (b) NSGAII, (b) BnB-NSGAII and (d) BnB-NSGAII legacy methods. Feasible and infeasible individuals are plotted in green and red respectively.}
\end{figure}

In this experiment, the evaluation of the performance of each method is limited to how many times the method finds at least 1 feasible solution over the 10 runs. Figure \ref{fig:stat} shows the number of times each method succeeded the test. It can be obviously concluded that BnB-NSGAII legacy method overcomes the performances of NSGAII and BnB-NSGAII. It should be noted that the computational effort is not regarded since all the runs converge within 30 minutes. Which is considered an acceptable time for such a problem.\par
Figure \ref{fig:nsga_failed} shows that NSGAII explored local space of the domain depending on the initial population. While Figure \ref{fig:bnb_nsga_failed} shows that BnB-NSGAII explored random spaces of the domain. Figure \ref{fig:bnb_legacy} shows that the legacy feature guides the exploration force of BnB-NSGAII towards the feasible solutions.\par
\section{Conclusion} \label{conc}
The 3 stage reducer problem is a point of interest of many researchers, either to use/ enhance it for engineering applications, or to examine the performance of optimization methods. The 3SR problem is desirable for such experiments for its complexity.\par
The 3SR problem was reformulated to a bi-objective problem in this paper to demonstrate a proposed application of BnB-NSGAII. The proposed application is to use BnB-NSGAII as an initiator of NSGAII, where BnB-NSGAII initially seeks feasible individuals before injecting them into the initial population of NSGAII. \par
BnB-NSGAII was enhanced by adding the legacy feature. The legacy feature is a generic feature that can be implemented in any branch and bound algorithm. Any parameter that is tuned during the node solving process could be the legacy. In this paper, the legacy was the last population in the father node in BnB-NSGAII. The latter was then used to initialize the child node.\par
The performances of NSGAII and BnB-NSGAII with and without the legacy feature were tested on the bi-objective version of the 3SR problem. Results show that the legacy feature guides the exploration force of BnB-NSGAII leading it to a better solution than that obtained by NSGAII and BnB-NSGAII.\par

\bibliographystyle{splncs04}
\bibliography{3stages_reducer}
\appendix
\section{3SR Problem Constraints} \label{prob_cnt}
\subsection{Closure condition}
Interference and fitting constraints are adopted from \cite{Han2020}. In \cite{Fauroux2004}, the closure condition was expressed with the distance between the terminal point $O_3$ shown in Figure \ref{fig:reducteuriftommvueface} and required position of the center of the output shaft. The coordinate of $O_3$ can be easily compute with the center distance of each stage and the angle $\xi_1$, $\xi_2$ and $\xi_3$. But, if we consider that center distance of each stage allow this closure condition, we can compute the value of $\xi_2$ and $\xi_3$. By this way can reduce he number of variables in the optimization problem.\par
For a given value of $\xi_1$ and $r_{1,1}$, $r_{1,2}$, center distance of each stage allow a closure if we have :
\[
\Vert\vec{O_1O_3}\Vert\leq\Vert\vec{O_1O_2}\Vert+\Vert\vec{O_2O_3}\Vert
\]
Assuming the previous condition is true, we can compute the two intersection of circle of center $O_1$ of $\Vert\vec{O_1O_2}\Vert$ radius and circle of center $O_3$ of $\Vert\vec{O_2O_3}\Vert$ radius.\par
With $a_2=\Vert\vec{O_1O_2}\Vert$ and $a_3=\Vert\vec{O_2O_3}\Vert$ we have :
\[
\left\lbrace
\begin{array}{l}
a_2\sin\alpha_1 - a_3\sin\alpha_3=0\\
a_2\cos\alpha_1 + a_3\cos\alpha_3=\Vert\vec{O_1O_3}\Vert
\end{array}
\right.
\]
which give :
\[
\cos\alpha_1 = \frac{\overline{O_1H}}{\overline{O_1O_2}}=\frac{a_2^2-a_3^2 + (O_1O_3)^2}{2a_2\overline{O_1O_3}}
\]

\begin{figure}[H]
\centering
\includegraphics[width=0.3\linewidth]{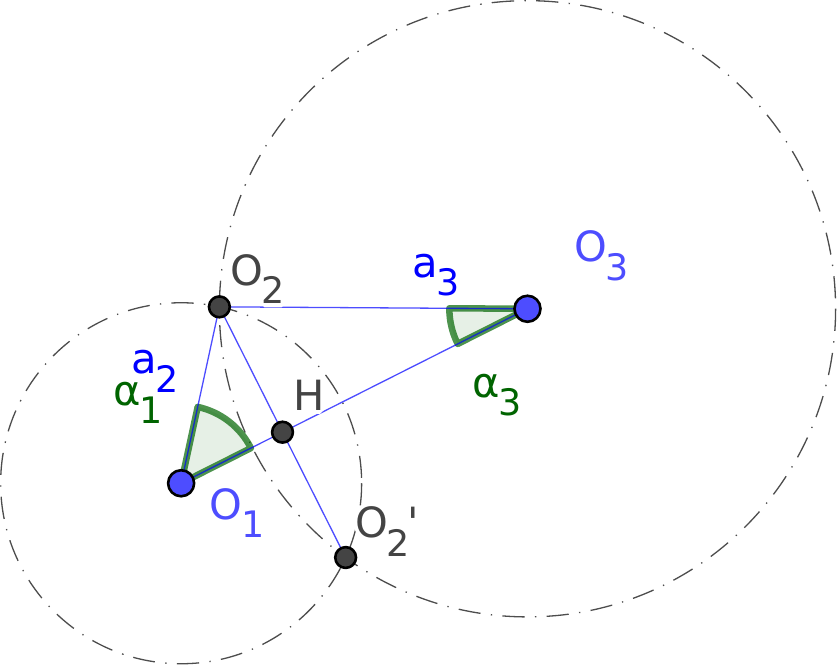}
\caption{}
\label{fig:closurefigure}
\end{figure}

Knowing $\alpha_1$, computation of coordinate of $O_2$ and $O_2^{'}$ is straightforward. if the two position $O_2$ and $O_2^{'}$ allow the wheel of the 2nd stage to fit in the casing box, then $O_2^{'}$ is preferred for lubrication reason.\par
\subsection{Mechanical constraint for one stage of the mechanism}
\subsubsection{Constraints related to the gear pair}
Following the recommendation from International Standard ISO 6336, \cite{ISO63361},\cite{ISO63362},\cite{ISO63363} we can calculate, knowing the geometry of gear pair, the material and the working conditions the contact and $\sigma_{\text{H}}$ the bending stress $\sigma_{\text{F}}$ in the gear pair. These stresses must be less of equal to the respective permissible value $\sigma_{\text{HP}}$ and $\sigma_{\text{FP}}$, depending on the material and the working conditions.\par
From \cite{ISO63363} the bending stress $\sigma_{\text{F}}$ is given by ($1$ for the pinion and $2$ for the wheel)::
\begin{equation*}
\sigma_{\text{F}(1,2)}=\sigma_{\text{F0}}\left(K_\text{A}K_\text{V}K_{\text{F}\alpha}K_{\text{F}\beta}\right)
\end{equation*}
with $\sigma_{\text{F0}(1,2)}$, the nominal tooth stress :
\begin{equation*}
\sigma_{\text{F0}(1,2)}=\frac{F_t}{bm_\text{n}}\left(Y_\text{F}Y_\text{S}Y_\beta Y_\text{B}Y_\text{DT}\right)
\end{equation*}
where :
\begin{itemize}
\item $F_t$ : is the tangential load from \cite{ISO63361}.
\item $b$ : is the facewidth.
\item $m_\text{n}$ : is the normal module.
\end{itemize}

Factors $K_\text{A}$, $K_\text{V}$, $K_{\text{F}\alpha}$, $K_{\text{F}\beta}$ are related to dynamic ad loading conditions in the gear. Factors $Y_\text{F}$, $Y_\text{S}$, $Y_\beta$, $Y_\text{B}$, $Y_\text{DT}$ are related to the geometry effect on stress.\par
From \cite{ISO63363}, the permissible bending stress $\sigma_{\text{FP}}$ is given by :
\begin{equation*}
\sigma_{\text{FP}}=\frac{\sigma_\text{FLim}}{S_\text{Fmin}}\left(Y_\text{ST}Y_\text{NT}Y_{\delta\text{relT}}Y_{R\text{relT}}Y_\text{X}\right)
\end{equation*}
with $\sigma_\text{FLim}$ is the nominal stress number (bending) from reference test gears \cite{ISO63365} and $S_\text{Fmin}$ the minimal required safety factor.
Factors $Y_\text{ST}$, $Y_\text{NT}$,  $Y_{\delta\text{relT}}$, $Y_{R\text{relT}}$, $Y_\text{X}$ are related to the reference test gears and the geometry and material conditions of the gear pair.\par
From \cite{ISO63362} the contact stress is given by ($1$ for the pinion and $2$ for the wheel):
\begin{equation*}
\sigma_{\text{H}(1,2)}=Z_{(\text{B,D})}\sigma_{\text{H0}}\sqrt{K_\text{A}K_\text{V}K_{\text{H}\alpha}K_{\text{H}\beta}}
\end{equation*}
with $\sigma_\text{H0}$ is the nominal contact stress :
\begin{equation*}
\sigma_{\text{H0}}=\left(Z_\text{H}Z_\text{E}Z_\varepsilon Z_\beta\right)\sqrt{\frac{F_t}{bd_1}\frac{u+1}{u}}
\end{equation*}
Factors $Z_\text{H}$, $Z_\text{E}$ ,$Z_\varepsilon$, $Z_\beta$ are related to the Hertzian theory of contact, and take into account geometry and material in the gear pair.\par
From \cite{ISO63362} the permissible contact stress $\sigma_{\text{HP}}$ is :
\begin{equation*}
\sigma_{\text{HP}}=\frac{\sigma_\text{HLim}}{S_\text{Hmin}}\left(Z_\text{NT}Z_\text{L}Z_\text{V}Z_\text{R}Z_\text{W}Z_\text{X}\right)
\end{equation*}
with $\sigma_\text{HLim}$ is the allowable contact stress number and $S_\text{Hmin}$ is the minimum required safety factor for surface durability.
Factors $Z_\text{NT}$, $Z_\text{L}$, $Z_\text{V}$, $Z_\text{R}$ ,$Z_\text{W}$, $Z_\text{X}$ are related to lubrication conditions, surface roughness and hardened conditions and size of the tooth.\par
So to respect the requirement specification of a given power to be transmitted, the gear pair must respect :
\begin{align*}
\sigma_{\text{F}(1,2)}&\leq\sigma_{\text{FP}}\\
\sigma_{\text{H}(1,2)}&\leq\sigma_{\text{HP}}
\end{align*}
Considering that $\sigma_{\text{F}}$ is proportional to $F_t$  and $\sigma_{\text{H}}$ is proportional to $\sqrt{F_t}$ for a given gear pair, we can rewrite these 2 conditions with $P_t$ the power to be transmitted :
\begin{align*}
\frac{\sigma_{\text{FP}}}{\sigma_{\text{F}(1,2)}}P_t&\geq P_t \\
\left(\frac{\sigma_{\text{HP}}}{\sigma_{\text{H}(1,2)}}\right)^2 P_t&\geq P_t
\end{align*}
Usually, some factors are slightly for the pinion and the wheel so transmitted power is different for the pinion $(1)$ and the wheel $(2)$. We will keep the minimal value.\par
So finally, for the stage number $s$ on the reducer, the following conditions must be fulfilled :
\begin{align}
\min\left(\frac{\sigma_{\text{FP}s}}{\sigma_{\text{F}(1,2)s}}\right)P_t&\geq P_t \\
\min\left(\frac{\sigma_{\text{HP}s}}{\sigma_{\text{H}(1,2)s}}\right)^2 P_t&\geq P_t
\end{align}
Following condition must be respected :
\begin{itemize}
\item For the transverse contact ratio : $\varepsilon_\alpha \geq 1.3$.
\item For the minimal face width : $b\geq 0.1d_2$
\item For the maximal face width : $b\leq d_1$
\end{itemize}

In order to use pinion with at least $Z_\text{min}=14$ teeth, the value of the profile shift coefficient must be adjusted to avoid gear meshing with the relation :
\[
Z_\text{min} \geq \frac{2(1-x_1)}{\sin\alpha_\text{n}^2} \Rightarrow x_1\geq 1-Z_\text{min}\frac{\sin\alpha_\text{n}^2}{2}\Rightarrow x_1\geq 0.1812
\]

\subsubsection{Constraint related to shaft's reducer}
In each of the 4 shafts of the mechanism, the transmitted torque produce shear stress. This stress must not exceed the allowable shear of the material of shafts $\tau_\text{max}$. We assume here that all the shaft are using the same steel and that all shaft can be consider as beam.
So, with $r_{a,0}$, the radius of input shaft, and $r_{a,s}$, $s=1 \ldots 3$ the radius of output shaft of the three stages, we have :
\begin{equation}
\tau_s=\frac{2C_s}{\pi {r_{a,s}}^3 }\leq \tau_\text{max}\text{ for }s=1 \ldots 3
\end{equation}
$C_s$ is the output torque of each stage and $C_e$ the torque on the input shaft, where $Z_{i,1}$ and $Z_{i,2}$ are the number of teeth for pinion $(_1)$ and wheel $(_2)$ of stage number $i$:
\[
C_s = C_e\prod\limits_{i=1}^{i=s}\frac{Z_{i,2}}{Z_{i,1}}
\]
For the input shaft we have :
\begin{equation}
\tau_0=\frac{2C_e}{\pi {r_{a,0}}^3 }\leq \tau_\text{max}
\end{equation}
The total rotation angle between the initial section of the input shaft and the final section of the output shaft is :
\[
\theta = \frac{2C_el_{a,0}}{G\pi {r_{a,0}}^3 } + \sum\limits_{s=1}^{s=3}\frac{2C_s l_{a,s}}{G\pi {r_{a,s}}^3 }
\]
For some reasons (dynamic behaviour of the reducer, ...) this total rotation angle should be limited by a maximal value $\theta_\text{max}$.
\begin{equation}
\theta \leq \theta_\text{max}
\end{equation}
\end{document}